\documentclass[reqno,11pt]{article}
\usepackage{times}
\usepackage[latin1]{inputenc}
\usepackage[T1]{fontenc}
\usepackage{amssymb}
\usepackage{amsmath}
\usepackage{amsthm}
\usepackage{latexsym}
\usepackage[dvips]{graphicx}
\usepackage[margin=2.5cm]{geometry}

\def\R{{\mathbb R}}

\def\e{{\varepsilon}}

\def\G{{\Gamma}}
\def\s{{\sigma}}

\def\d{{\delta}}

\def\p{{\prime}}

\def\F{ {\mathcal{F} } }

\def\S{{Schr\"{o}dinger }}

\def\pr{\vskip3pt \nl{\bf Proof. \,\,}}

\def\+R{+_{_{ \!\! \R}}}

\DeclareMathAlphabet{\mathpzc}{OT1}{pzc}{m}{it}
\def\nl{\vglue0.3truemm\noindent}
\numberwithin{equation}{section}

\begin{document}
\newtheorem{theo}{Theorem}[section]
\newtheorem{pro}[theo]{Proposition}
\newtheorem{lem}[theo]{Lemma}
\newtheorem{defin}[theo]{Definition}
\newtheorem{rem}[theo]{Remark}
\newtheorem{cor}[theo]{Corollary}

\title{A new proof of long range scattering for \\ critical nonlinear \S equations}

\author{ 
Jun Kato\footnote{JSPS Postdoctoral Fellow for Research Abroad.}
\\
\vspace{-.2truecm}
{\footnotesize Graduate School of Matemathics}
\\{\footnotesize Nagoya University}
\vspace{-.2truecm}
\\{\footnotesize Furocho, Chikusaku, Nagoya (Japan) }
\vspace{-.2truecm}
\\{\scriptsize jkato@math.nagoya-u.ac.jp}
\\
\and
Fabio Pusateri\\
\vspace{-.2truecm}
{\footnotesize Courant Institute of Mathematical Sciences }
\\{\footnotesize New York University}
\vspace{-.2truecm}
\\{\footnotesize 251 Mercer Street, N.Y., N.Y. 10012  (USA) }
\vspace{-.2truecm}
\\{\scriptsize pusateri@cims.nyu.edu}
}

\date{}

\maketitle

\begin{abstract}
We present a new proof of global existence and long range scattering, from small initial data, 
for the one-dimensional cubic gauge invariant nonlinear Schr\"odinger equation,
and for Hartree equations in dimension $n \geq 2$.
The proof relies on an analysis in Fourier space, related to
the recent works of Germain, Masmoudi and Shatah
on space-time resonances.
An interesting feature of our approach is that we are able to identify
the long range phase correction term
through a very natural stationary phase argument.
\end{abstract}

\section{Introduction}

The problem of asymptotic behavior of small solutions of nonlinear evolution equations
has been extensively treated by many authors in the past fourty years. 
Almost all of the entire literature on the subject is dedicated to 
prove that solutions of nonlinear PDEs, evolving from small Cauchy data, 
behave asymptotically in time like solutions of the associated linear problem. 
The task of identifying nonlinear global dynamics is extremely challenging.
Indeed, there are not many results in the field where a small nonlinear solution
is proven to exist globally and behave, for large times, differently than a linear one.

Among the most celebrated results concerning nonlinear \S equations we want to mention 
the works of Deift and Zhou\footnote{
The results of these authors, based on the scattering/inverse scattering theory,
also apply to large data.
} (see \cite{DZ} and references therein)
and the work of Hayashi and Naumkin \cite{HN}.
These latter authors also proved long range scattering results for other dispersive equations
such as Benjamin-Ono \cite{HNBO} and KdV \cite{HNKdV}.
For some results concerning nonlinear asymptotic behavior of solution of wave equations see
\cite{LR1,Alinhac1,KK} and references therein.

All of the methods employed to deal with the problems mentioned above
do not suggest any unified approach to the question of asymptotic behavior
of small solutions when ordinary (linear) scattering fails.
Our aim, in this short note, is to present a first application of a new simple method
that can be used in situations where the asymptotic behavior of a nonlinear solution
differs from the linear one by a phase correction term.

More specifically, we are going to show a new proof of the results contained in \cite{HN}
about long range scattering 
for the one dimensional cubic gauge-invariant NLS equation,
and for Hartree equations in dimension $n \geq 2$.
These equations provide us with an easy setting to expose and implement our method
since the necessary computations turn out to be particularly straightdforward.
Nevertheless, we believe that our approach is fairly general,
and could be used in similar situations.
In particular, we think it could be applicable to other scattering critical equations
where modified wave operators can be -- or have been -- constructed with the help 
of a phase correction term,
but asymptotics for solutions of the Cauchy problem are not known yet.

As we will explain below, there are substantial differences 
between our method and the one used in \cite{HN}.
In particular, we believe that the main contributions of our approach lie 
in its simplicity and in the way we explicitely derive the phase correction term.
No a priori knowledge of the latter is required,
and its derivation is just a mere consequence of a very natural application
of the classical stationary phase lemma.

The setting for our proof is inspired from the recent works 
of Germain, Masmoudi, and Shatah \cite{GMS1, GMS2, GMS3} on the analysis of space-time resonances.
The very general method developed by these authors,
combining earlier ideas about normal forms and commuting vector fields
with tools from harmonic analysis,
has already proven to be extremely useful in the study of global existence 
of asymptotically free solutions to nonlinear dispersive equations.
In this note we show how a similar approach 
can be also be applied in the context of long range scattering.

\vskip5pt
Before formulating the problems and explaining the main ideas of the proof,
we give the following
\noindent
{\bf Notations}:
For $g \in \mathcal{S}(\R^{n})$, we denote its Fourier transform by 
\begin{equation}
\label{FT}
  \widehat{g}(\xi)=\mathcal{F}[g](\xi)
  = {(2\pi)}^{-n/2}\int e^{-ix\cdot\xi}g(x)\,dx \, .
\end{equation}
To avoid ambiguities we will sometime use subscripts to indicate the variables with respect to which 
the Fourier transform is taken.
For $m,l\in\R$ we define weighted Sobolev spaces $H^{m,l}$ and $\dot{H}^{m,l}$ by 
\begin{align} 
\label{Hms}
  H^{m,l} & = \left\{ \varphi\in\mathcal{S}' \, :  {\|\varphi\|}_{H^{m,l}}:={\|\langle x\rangle^{l}\langle \nabla	
  \rangle^{m}\varphi\|}_{L^{2}}<\infty \right\}  \, ,
\\
\nonumber
  \dot{H}^{m,l} & = \left\{ \varphi\in\mathcal{S}' \, :  
  { \| {| x |}^{l} {| \nabla |}^{m} \varphi\|}_{L^{2}}<\infty \right\} 
\end{align}
where $\langle x \rangle := {(1+{|x|}^{2})}^{1/2}$, 
$\langle \nabla \rangle:=\mathcal{F}^{-1}\langle \xi \rangle\mathcal{F}$
and 
$| \nabla | := \mathcal{F}^{-1} | \xi | \mathcal{F}$.
\noindent
We denote by $e^{i t \Delta/2}$ the fundamental solution of the linear \S equation $i u_t + \frac{1}{2} \Delta u = 0$
and define
\begin{equation}
\label{profile}
f (t,x) := e^{-it\Delta/2} u(t,x) 
\end{equation}
to be the {\it profile} of a solution $u$ of \eqref{eqs}.

\subsection{The problems}
We consider the Cauchy problem for the nonlinear \S and Hartree equations
\begin{equation}
\label{eqs}
\left\{
\begin{aligned}
 & i \partial_{t}u + \frac{1}{2}\Delta u
   = g(u) \quad, \quad  (t,x) \in \R \times \R^n \quad,
   \\
 & u|_{t=1} = e^{i \Delta/2} u_{\ast} 	\quad  , \quad  x \in \R^{n} \quad,
\end{aligned}
\right.
\end{equation}
where the nonlinear term is 
\begin{equation}
\label{gNLS}
  g(u) =  {|u|}^{2} u  \quad, \quad x \in \R
\end{equation}
in the NLS case, and
\begin{equation}
\label{gH}
  g(u) =  \left( {|x|}^{-1}\ast {|u|}^2 \right) u \quad , \quad  x \in \R^{n} \quad , \quad n \geq 2
\end{equation}
in the Hartree case.
The nonlinearities \eqref{gNLS} and \eqref{gH} are critical
from the point of view of long time asymptotic behavior,
since the $L^2$--norm of the nonlinear terms, computed on a linear solution, decays like $t^{-1}$. 
Moreover, they belong to the case of long range scattering
where ordinary wave operators fail to exist \cite{Barab,HT}.
There are many papers treating the problem of long time behaviour for
solutions of \eqref{eqs} with \eqref{gNLS} or \eqref{gH}.
Modified wave operators including a certain phase factor
were first constructed by Ozawa \cite{Ozawa} in the case of \eqref{gNLS}.
This result was then extended to other critical NLS equations and Hartree equations 
by Ginibre and Ozawa \cite{GO}.
Asymptotic completness for \eqref{eqs}--\eqref{gNLS} was shown by Naumkin \cite{Naumkin},
for small data in $H^{1,0} \cap H^{0,1}$. 
Subsequently, Hayashy and Naumkin \cite{HN} proved asymptotic completness for 
scattering critical NLS equations in dimensions up to $3$, and for Hartree equations in dimension $n \geq 2$.
They also provided a precise asymptotic formula for solutions of \eqref{eqs} with \eqref{gNLS} or \eqref{gH};
see, respectively, \eqref{NLSas} and \eqref{Has} below.

\subsection{Main ideas in the proof}
When trying to obtain global existence of small solutions
to scattering critical nonlinear dispersive equations,
the main issue one often has to face
is to establish the sharp (linear) decay of the $L^{\infty}_x$ norm of solutions.
Before the work of Hayashi and Naumkin \cite{HN}, 
no author was able to prove that solutions of \eqref{eqs} with \eqref{gNLS},
respectively \eqref{gH},
satisfy the linear decay estimate $\|u(t)\|_{L^\infty} \lesssim t^{-\frac{1}{2}}$,
respectively $\|u(t)\|_{L^\infty} \lesssim t^{-\frac{n}{2}}$.
In order to obtain these bounds, Hayashi and Naumkin introduced an amplitude modulation term 
-- commonly referred to as a ``phase correction'' --
which enabled them to show the desired $L^\infty$ decay.
More specifically, in \cite{HN} the authors consider the profile $f$ of a solution $u$, 
defined as \eqref{profile}, 
and write
\begin{align}
\nonumber
\partial_t \widehat{f} = e^{it \frac{{|\xi|}^2}{2}} \F [ g(u) ]
\\
\label{}
=  t^{-1}  g ( \widehat{f} ) + R
\label{HNdthatf}
\end{align}
where
$
R = e^{it \frac{{|\xi|}^2}{2}} \F [ g(u) ] - t^{-1} g ( \widehat{f} )
$.
Subsequently, they perform some manipulations, based on a particular identity for the linear \S semigroup,
in order to show that $R$ decays faster than $t^{-1}$.
The first term in \eqref{HNdthatf} is then eliminated
by considering a modified profile $\widehat{w} (t,\xi) = \widehat{f} (t,\xi) e^{i B(t,\xi)}$,
for a suitable real valued phase function $B$.
The resulting equation yields an a priori $L^\infty$ bound on $\widehat{f}$,
which in turn implies the linear decay of solutions.

To understand the difference with our approach,
let us briefly describe how we are going to derive the phase correction term.
For the sake of discussion we just consider the case of \eqref{eqs} with \eqref{gNLS}.
We refer to section \ref{secproNLS} for a complete proof in this case
and to section \ref{secproH} for the proof in the Hartee case.
Following the already cited works of Germain, Masmoudi, and Shatah,
we write Duhamel's formula in Fourier space for the profile $f$:
\begin{align}
\label{int0}
\widehat{f}(t,\xi) & = \widehat{u}_{\ast}(\xi)
    - \frac{i}{2\pi} 
    \int_{1}^{t}  \int e^{is \eta (\xi-\sigma)}
    \widehat{f}(s,\xi-\eta) \overline{\widehat{f}(s,\sigma-\eta)}\widehat{f}(s,\sigma)\,d\eta d\sigma ds  \, .
\end{align}
We then apply a stationary phase argument
to the oscillatory integral 
with respect to the variables $\eta$ and $\sigma$, 
obtaining
\begin{align}
\nonumber
\widehat{f} (t, \xi) =  \widehat{u}_{\ast}(\xi)
- i  \int_1^t \frac{1}{s} \,
{ |  \hat{f} (s, \xi) | }^2  \hat{f} (s, \xi) \, ds  + \int_1^t R (s, \xi) \, ds \, ,
\end{align}
which is an integral version of \eqref{HNdthatf}.
Since $R$ can be easily shown to decay faster that $s^{-1}$,
we have explicitely determined, in a natural and direct fashion, 
the phase correction term which is necessary in the construction of modified wave operators.

Another interesting feature of our approach
lies in the way we estimate weighted norms of $f$.
We completely avoid the use of the vector field\footnote{
We use \eqref{int1} and \eqref{int3} 
to estimate directly $\partial^k_\xi \widehat{f}$ in $L^2$, which correspond to estimating $J^k u$ in $L^2$.
} 
$J = x - i t \nabla$,
and the identity $J = e^{it \frac{\Delta}{2}} x e^{-it \frac{\Delta}{2}}$.
This vector field approach, despite being customary in the literature, gives, in our opinion,
a more restricted perspective on the problem.
The Fourier side approach 
provides, instead, a more general framework, 
where the particular structure of the gauge invariant non-linearity manifests itself naturally.
To see this, let us describe how we establish weighted $L^2$ bounds on $f$.
Since weights in physical space correspond to derivatives on the Fourier side,
we apply $\partial_\xi$ to the integral equation \eqref{int0}
associated to \eqref{eqs}--\eqref{gNLS}. 
Then, the main issue is to estimate the contribution of the term where $\partial_\xi$ hits the 
oscillating phase $e^{i s \varphi}$, with $\varphi =  \eta (\xi-\sigma)$,
because this results in the appearence of a growing factor of $s$.
It is at this stage that the structure of the gauge invariant nonlinearity comes into play.
Observe, indeed, that differentiating $e^{i s \varphi}$ also introduces 
a factor $\partial_\xi \varphi = \eta$
that vanishes exactly on the set where $\partial_\s \varphi$ vanishes.
This suggests that, integrating by parts with respect to $\s$, 
one should be able to recover a factor of $s^{-1}$.
In this particularly simple case, one just notices that
$\partial_\xi e^{is\eta\cdot(\xi-\sigma)} =  - \partial_\s e^{is\eta\cdot(\xi-\sigma)}$,
so that, after integration by parts, no extra--growth in $s$ is actually present.
Notice that the same type of argument would not work for the other cubic nonlineairties $u^3$, $\bar{u}^3$
or $\bar{u}^2 u$.
Although the present example is quite elementary,
this type of analysis of ``generalized null-structures''
constitutes an important part of the methods developed by Germain, Masmoudi and Shatah.

\subsection{Reformulation of \eqref{eqs} with \eqref{gNLS} and statement of results}
The integral equation associated to \eqref{eqs}--\eqref{gNLS} is
\begin{align}
\nonumber
u(t) & =  e^{i (t - 1 ) \partial_{xx}/2 }  u(1)
    -i  \int_{1}^{t} e^{i(t-s) \partial_{xx}/2}
    {|u(s)|}^{2} u(s) \,ds  
\\
\label{intu}
     & =:  e^{i t \partial_{xx}/2 } u_{\ast}  +  e^{it \partial_{xx}/2 } C(u,u,u)(t) \, .
\end{align}
If $f$ denotes the profile of $u$ defined in \eqref{profile},
then
\begin{equation*}
  f(t) =  u_{\ast}
    -i  \int_{1}^{t}   e^{- is \partial_{xx}/2}  |u(s)|^{2} u(s) \,ds \, .
\end{equation*}
In Fourier space the above equation becomes\footnote{
With the normalization in \eqref{FT}
one has $ \F (f \ast g) = {(2\pi)}^{n/2} \F (f) \F (g)$
and
$ \F (f g) = {(2\pi)}^{-n/2} \F (f) \ast \F (g)$.
}
\begin{align}
\nonumber
\widehat{f}(t,\xi) & = \widehat{u}_{\ast}(\xi)
    -i  {(2\pi)}^{-1} 
    \int_{1}^{t}  \int e^{is \eta (\xi-\sigma)}
    \widehat{f}(s,\xi-\eta) \overline{\widehat{f}(s,\sigma-\eta)}\widehat{f}(s,\sigma)\,d\eta d\sigma ds  
\\
\label{int}
    & =:  \widehat{u}_{\ast}(\xi)  +  \widehat{C} (u,u,u)(t)
    \, ,
\end{align}
since
\begin{equation*}
\frac{1}{2}  \left[ \xi^{2} - (\xi-\eta)^{2} + (\sigma-\eta)^{2} - \sigma^{2} \right] =  \eta (\xi - \sigma) \, .
\end{equation*}
Changing variable in the $\s$--integral, \eqref{int} can also be written as
\begin{equation}
\label{int1}
  \widehat{f}(t,\xi) =  \widehat{u}_{\ast}(\xi)
    -i  {(2\pi)}^{-1} 
    \int_{1}^{t}  \int e^{i s \eta  \sigma}
    \widehat{f}(s,\xi-\eta) \overline{\widehat{f}(s,\xi - \sigma-\eta)}
    \widehat{f}(s,\xi - \sigma)\,d\eta d\sigma ds  \, .
\end{equation}

\nl
We are going to prove the following:
\begin{theo}
\label{theoNLS}
Let $u_\ast \in H^{1,0} \cap H^{0,1}$ with ${\| u_\ast \|}_{H^{1,0}} + {\| u_\ast \|}_{H^{0,1}} \leq \e$ 
for some $\e$ small enough.
Then there exists a unique global solution 
$u \in C (\R ; H^{1,0}(\R) \cap H^{0,1}(\R))$ of \eqref{eqs}--\eqref{gNLS} such that
\begin{equation}
\label{decayNLS}
{\|  u (t) \|}_{L^\infty} \lesssim \e {(1 + |t|)}^{-1/2} \, .
\end{equation}
Furthermore, there exists a unique function
$W \in L^\infty$
such that for $t \geq 1$
\begin{equation}
\label{scatteringNLS}
{\left\|  \widehat{f}  \exp \left( - i \int_1^t  {|\widehat{u}(s)|}^2 \, \frac{ds}{s} \right)  
																	- W  \right\|}_{L^\infty} \lesssim  t^{-\d}
\end{equation}
for some $0 < \d < \frac{1}{2}$, where $f$ is defined by \eqref{profile}.
\end{theo}

\nl
{\bf Remarks}.
\begin{enumerate}
\item
Theorem \ref{theoNLS} corresponds to the $n=1$ case of theorem 1.1 and (the first part of) theorem 1.2 in \cite{HN}.

\item
It is not hard to see that the function $W$ in theorem \ref{theoNLS} also belongs to $L^2$.
As a consequence, \eqref{scatteringNLS} also holds with the $L^2$ norm replacing the $L^\infty$ norm. 
We will not provide the details here, since it will be clear
how to derive this fact from our computations.
\item
The estimate \eqref{scatteringNLS} 
is a direct byproduct of the proof of the $L^\infty$ decay estimate \eqref{decayNLS}.
Furthermore, starting from \eqref{scatteringNLS},
the same arguments of \cite[381-383]{HN}
can be used to establish the following asymptotic behaviour of solutions of \eqref{eqs}--\eqref{gNLS}:
there exists a unique function $\Phi \in L^\infty$ such that as $t \rightarrow \infty$
\begin{equation}
\label{NLSas}
u (t,x) = \frac{1}{ {(i t)}^{1/2} } W ( x/t ) 
											\exp \left( i \frac{{|x|}^2}{4t} + i {| W(x/t) |}^2 \log t + i \Phi (x/t) \right)
											+ O \left( t^{- 1/2  - \beta} \right) \, 
\end{equation}
uniformly in $x \in \R$, for some small $\beta > 0$.
\end{enumerate}

\nl
In analogy with \cite{HN}, we define the space
\begin{equation*}
X_T := \left\{ u \, :  {\| u \|}_{X_T} :=  {\| t^{1/2} u \|}_{L^\infty_T L^\infty_x}  + 
																		 {\| t^{-\alpha} u \|}_{L^\infty_T \dot{H}^{1,0}_x} 
																	+  {\| t^{-\alpha} f \|}_{L^\infty_T \dot{H}^{0,1}_x} 
																	+  {\| u \|}_{L^\infty_T L^2_{x}}  < \infty     \right\} 
\end{equation*}
where $L^\infty_T := L^\infty_t ([0,T])$, $f$ is defined in \eqref{profile} 
and $\alpha > 0$ is small enough.
The global solution of theorem \ref{theoNLS} will be constructed 
as a continuation of the local solution provided by the following classical result:

\begin{theo}[Local existence \cite{CW,GV1,HT2}]
\label{theolocal}
Assume ${\| u_\ast \|}_{H^{1,0}} + {\| u_\ast \|}_{H^{0,1}} \leq \e$, for some  $\e$ small enough.
Then there exists a finite time $T > 1$ and a unique solution 
$u \in C ([0,T]; H^{0,1} (\R) \cap H^{1,0} (\R))$ of \eqref{eqs}--\eqref{gNLS},
such that ${\| u \|}_{X_T} \leq C \e$.
\end{theo}

\nl
The first part of theorem \ref{theoNLS} will then result as a consequence of
\begin{pro}
\label{proNLS}
Let $u$ be given as in theorem \ref{theolocal}, 
then
\begin{equation}
\label{cubicest}
{\left\| u \right\|}_{X_T} \leq \e + C {\| u \|}^3_{X_T}  \, ,
\end{equation}
for some constant $C$ independent of $T$.
\end{pro}
\nl
The proof of proposition \ref{proNLS} is performed in section \ref{secproNLS}.

\subsection{Reformulation of \eqref{eqs} with \eqref{gH} and statement of results}
Let us consider now the Cauchy problem for the Hartree equation with a Coulomb potential
\eqref{eqs}--\eqref{gH}.
The corresponding integral equation is
\begin{align}
\nonumber  
		u(t) & = e^{i(t-1) \Delta/2} u (1)
    -i \int_{1}^{\infty} e^{i(t-s)\Delta/2}
    \left(  {|x|}^{-1} \ast {|u(s)|}^{2} \right) u(s) \, ds 
    \\
    \label{CH}
    & =:  e^{i t \Delta/2 } u_{\ast}  +  e^{it \Delta/2 } C(u,u,u)(t)
\end{align}
so that the profile $f$ satisfies
\begin{align}
\label{intH}
\widehat{f}(t,\xi) & = \widehat{u}_{\ast}(\xi)
    -i {(2\pi)}^{-n/2}
    C_1  \int_{1}^{t}\int e^{is\eta\cdot(\xi-\sigma)}
    {|\eta|}^{-n+1} \widehat{f}(s,\xi-\eta) 
    \overline{\widehat{f}(s,\sigma-\eta)}\widehat{f}(s,\sigma) \, d\eta d\sigma ds 
\\
\nonumber
& = \widehat{u}_{\ast}(\xi) + \widehat{C} (u,u,u) \, ,
\end{align}
since
$$
\frac{1}{2} \left[ {|\xi|}^{2}-{|\xi-\eta|}^{2} + {|\sigma-\eta|}^{2} - {|\sigma|}^{2} \right] 
= \xi\cdot\eta-\sigma\cdot\eta \, .
$$
Here\footnote{
We are using the fact that
\begin{equation*}
\F \left[ {|x|}^{-a} \right] = {2}^{n/2-a} \frac{ \Gamma((n-a)/2) }{ \G(a/2) }  {|x|}^{a - n}
									=: C_a {|x|}^{a-n}  
									\, .
\end{equation*}
Notice that $C_1 C_{n-1} = 1$.
} 
$C_1 := {2}^{n/2-1} \pi^{-1/2} \Gamma((n-1)/2)$.
After a change of variables, the integral equation 
\eqref{intH} can be written as
\begin{equation}
\label{int3}
  \widehat{f}(t,\xi) = \widehat{u}_{\ast}(\xi)
    - i {(2\pi)}^{-n/2} C_1 \int_1^t \int e^{is\eta\cdot\sigma}  {|\eta|}^{-n+1} \widehat{f}(s,\xi-\eta) 
    \overline{\widehat{f}(s,\xi - \sigma - \eta)}\widehat{f}(s,\xi - \sigma) 
    \,d\eta  d\sigma ds \, . 
\end{equation}

\nl
We are going to prove the following:
\begin{theo}
\label{theoH}
Let $m := [n/2] + 1$.
Assume that $u_\ast \in H^{m,0} \cap H^{0,m}$ and ${\| u_\ast \|}_{H^{m,0}} + {\| u_\ast \|}_{H^{0,m}} \leq \e$ 
for $\e$ small enough.
Then there exists a unique global solution $u \in C (\R ; H^{m,0}(\R^n) \cap H^{0,m}(\R^n) )$ 
of \eqref{eqs}--\eqref{gH}, such that
\begin{equation}
\label{decayH}
{\|  u (t) \|}_{L^\infty} \lesssim \e {(1 + |t|)}^{-n/2} \, .
\end{equation}
Furthermore, there exists a unique function $W \in L^\infty$
such that for $t \geq 1$
\begin{equation}
\label{scatteringH}
{\left\|  \widehat{f}  \exp \left( - i \int_1^t  \left( {|x|}^{-1} \ast {| \widehat{u} |}^2 \right) (s) \, \frac{ds}{s} \right)  
																	- W  \right\|}_{L^\infty} \lesssim  t^{-\d}
\end{equation}
for some $0 < \d < \frac{1}{2}$.
\end{theo}

\nl
Theorem \ref{theoH} above corresponds to Theorem 1.1$^\p$ and (part of) Theorem 1.2$^\p$ in \cite{HN}.
The remarks made after theorem \ref{theoNLS} 
apply also to the function $W$ in theorem \ref{theoH} and to the inequality \eqref{scatteringH}.
Furthermore, the same arguments in \cite[387-388]{HN}
can again be used to determine the asymptotic behaviour of solution of \eqref{eqs}--\eqref{gH}, i.e.
as $t \rightarrow \infty$
\begin{equation}
\label{Has}
u (t,x) = \frac{1}{ {(i t)}^{n/2} } W ( x/t ) 
											\exp \left( i \frac{{|x|}^2}{4t} + i \left( {|x|}^{-1} \ast {| W |}^2 \right) (x/t)  \log t 
											+ i \Phi (x/t) \right)
											+ O \left( t^{- n/2  - \alpha} \right) \,
\end{equation}
for some uniquely determined function $\Phi \in L^\infty$,
uniformly for $x \in \R^n$.

\nl
For the proof of theorem \ref{theoH}
we are going to follow the same strategy adopted in the proof of theorem \ref{theoNLS}.
Again in analogy with \cite{HN} we define
\begin{equation*}
X_T := \left\{  u \, : \, {\| u \|}_{X_T} := {\| t^{n/2} u \|}_{L^\infty_T L^\infty_x}
    + {\| t^{-\alpha}  u\|}_{L^\infty_T \dot{H}^{m,0}_x}  
    + {\| t^{-\alpha} f \|}_{L^\infty_T \dot{H}^{0,m}_x}
    + {\| u\| }_{L^\infty_T L^2_x}  
    < \infty
    \right\} \, .
\end{equation*}

\nl
Global solutions to \eqref{eqs}--\eqref{gH} will be constructed as a continuation of local solutions provided by

\begin{theo}[Local existence \cite{GV2,HT,HT2}]
\label{theolocalH}
Assume ${\| u_\ast \|}_{H^{m,0}} + {\| u_\ast \|}_{H^{0,m}} \leq \e$, for some  $\e$ small enough.
Then there exists a finite time $T > 1$ and a unique solution 
$u \in C ([0,T]; H^{0,m} (\R^n) \cap H^{m,0} (\R^n))$ of \eqref{eqs}--\eqref{gH},
such that ${\| u \|}_{X_T} \leq C \e$.
\end{theo}

\nl
In order to prove the first part of Theorem \ref{theoH}, 
it will then be sufficient to establish the following a priori estimates:
\begin{pro}
\label{proH}
Let $u$ be given as in \ref{theolocalH}. Then, there exists a constant $C$, independent of $T$, such that
\begin{equation}
\label{cubicestH}
{\left\| u \right\| }_{X_T} \leq \e + C {\| u \|}^3_{X_T}  \, .
\end{equation}
\end{pro}
\nl
The proof of proposition \ref{proH} is performed in section \ref{secproH}.

Before turning to the proofs, we recall some standard properties of the linear \S semigroup.
\begin{lem}[\cite{HN}]
\label{statNLS}
The \S semigroup satisfies the identity
\begin{equation*}
( e^{ it \Delta/2 } g  ) (x)  =  \frac{1}{ {(it)}^{n/2} } 
																				e^{i \frac{{|x|}^2}{2t} } \widehat{g} \left( x/t \right)
																				+ \frac{1}{t^{n/2 + \beta}} O \left( {\|  g \|}_{ H^{0, \gamma} } \right)
\end{equation*}
for $x \in \R^n$ and any $\gamma > \frac{n}{2}  + 2 \beta$.
In particular, it follows 
\begin{equation*}
{\| e^{ i t \Delta/2 } g  \|}_{L^\infty}  \lesssim  \frac{1}{ t^{n/2} } 
																				{\| \widehat{g} \|}_{L^\infty}
																				+ \frac{1}{t^{n/2 + \beta}} {\|  g \|}_{ H^{0, \gamma} }\, .
\end{equation*}
Furthermore, for any $ 2 \leq p \leq \infty$ one has the dispersive estimate:
\begin{equation}
\label{dispest}
{\| e^{ i t \Delta/2 } g \|}_{L^p}  \lesssim  t^{ -n (1/2 - 1/p) }  {\| g \|}_{L^{p^\p}}  \, .
\end{equation}
\end{lem}

\nl
{\bf Remark}.
Interpolating between the $L^2_x$ and $L^\infty_x$ components of the $X_T$ norm,
we have\footnote{
$L^{p,q}$ denotes the usual Lorentz space \cite{BL}.
}
\begin{equation}
\label{Lorentzdecay}
{\| u \|}_{L^{p,q}} \lesssim  t^{ -n (1/2 - 1/p) } {\| u \|}_{X_T} \, ,
\end{equation}
for any $2 < p < \infty$ and $1 \leq q \leq \infty$.

\vskip10pt
{\bf Acknowledgements.}
The authors are indebted with Prof. Jalal Shatah for the many enlightening discussions on the topic.
The first author is also grateful to the Courant Institute of Mathematical Science for the kind hospitality.


\section{Proof of proposition \ref{proNLS}}
\label{secproNLS}

Notice first that a solution to \eqref{eqs} with \eqref{gNLS} or \eqref{gH} enjoys conservation of the $L^2$ norm. 

\nl
{\it $\dot{H}^{1,0} \cap \dot{H}^{0,1}$ estimates}: 
From \eqref{intu} we have
\begin{eqnarray*}
{\| \partial_x C \|}_{L^2}   \lesssim \int_1^t  {\|\partial_x  u (s)\|}_{L^2}  {\|u(s)\|}^2_{L^\infty}  \, ds
														\lesssim    {\|u\|}^3_{X_T}  \int_1^t  s^\alpha s^{-1}  \, ds \lesssim {\|u\|}^3_{X_T} t^\alpha
\end{eqnarray*}
To estimate $ {\| x  C \|}_{L^2} = { \|\partial_\xi  \widehat{C} \|}_{L^2} $
we use \eqref{int1}. Then
\begin{align*}
&  \partial_\xi  \int_{1}^{t}  \int e^{is \eta \sigma}
    \widehat{f}(s,\xi-\eta) \overline{\widehat{f}(s,\xi - \sigma - \eta)} \widehat{f}(s,\xi - \sigma)
    \,d\eta d\sigma \, ds     
    \\
    & =  
    \int_{1}^{t}  \int e^{is \eta \sigma}
    \partial_\xi \widehat{f}(s,\xi-\eta) 
    	\overline{\widehat{f}(s,\xi- \sigma-\eta)} \widehat{f}(s,\xi - \sigma) \,d\eta d\sigma \,ds     
    	\hskip6pt + \hskip6pt \mbox{similar terms}
\end{align*}
where ``similar terms'' denotes the terms where the derivative with respect to $\xi$ falls on
the Fourier transform of the other profiles. These terms clearly enjoy the same estimates as the first one above.
Therefore, redistributing the phases on the three profiles, we see that
\begin{align*}
{ \| \partial_\xi \widehat{C} \|}_{L^2} & =  {\left\| \int_1^t e^{- is\partial_{xx}} 
\left( e^{is\partial_{xx}/2} (x f)  e^{-is\partial_{xx}/2} \bar{f}  e^{is\partial_{xx}/2} f  \right)  \, ds \right\|}_{L_2}
\\
& \lesssim  \int_1^t {\|u\|}^2_{L^\infty} {\|x f\|}_{L^2} \, ds 
\\
& \lesssim  \int_1^t  s^{-1} s^\alpha \, ds  \, {\|u\|}^3_{X_T}
\lesssim {\| u \|}_{X_T}^3 t^\alpha \, .
\end{align*}

\vskip5pt
\nl
{\it $L^\infty$ estimate}:
As we already pointed out, the advantage of writing \eqref{eqs}--\eqref{gNLS}
as an integral equation for the profile, 
is that one gets a clear understanding 
of why an asymptotic behavior such as \eqref{scatteringNLS} -- and eventually as \eqref{NLSas} -- occurs.
Indeed, following the proof of the classical stationary phase lemma, 
we use Plancharel's identity in \eqref{int1} to obtain\footnote{
We use the identity
$ \F \left[ \exp\left( i \langle x, Q x \rangle /2 \right)   \right] =
	e^{i \frac{\pi}{4} \mbox{sign } Q} {|\det{Q}|}^{-1/2} \exp\left( -i \langle \xi, Q^{-1} \xi \rangle/2 \right)
$,
with $Q =  \left( \begin{array}{cc} 0 & s \\ s & 0 \end{array} \right)$.
}
\begin{align*}
\widehat{f} (t, \xi)  & =  \widehat{u}_{\ast}(\xi)
    										- i   {(2\pi)}^{-1}  
    										\int_{1}^{t}  \int \F_{\eta, \sigma} \left[ e^{i s \eta \sigma } \right]
    										\F^{-1}_{\eta, \sigma} [F] (s, \eta, \sigma, \xi)
    										\, d\eta d\sigma ds
												\\
											& =  \widehat{u}_{\ast}(\xi)
    										- i  {(2\pi)}^{-1}  
    										\int_{1}^{t}  \int \frac{1}{s} \, e^{- i \eta \sigma/s }
    										\F^{-1}_{\eta, \sigma} [F] (s, \eta, \sigma, \xi)
    										\, d\eta d\sigma ds
\end{align*}
where
\begin{equation}
\label{F}
F(s, \eta, \sigma, \xi)  := \widehat{f}(s,\xi-\eta) \overline{\widehat{f}(s,\xi - \sigma-\eta)}
    												\widehat{f}(s,\xi - \sigma)  \, .
\end{equation}
Then
\begin{align}
\nonumber
\widehat{f} (t, \xi)  =   \widehat{u}_{\ast}(\xi)
    										& - i  {(2\pi)}^{-1}  \int_{1}^{t}  \frac{1}{s}  \int
    										\F_{\eta, \sigma}^{-1} [F] (s, \eta, \sigma, \xi) 
    										\, d\eta d\sigma ds
    										\\
    										\nonumber
    										&  - i  {(2\pi)}^{-1}  \int_{1}^{t}  \frac{1}{s} \int \left( e^{-i \eta \sigma/s } - 1 \right) 
    										\F_{\eta, \sigma}^{-1} [F] (s, \eta, \sigma, \xi)
    										\, d\eta d\sigma ds
    										\\
    										\label{intidfhat}
    										= 
    										\widehat{u}_{\ast}(\xi)
    										& - i  \int_1^t \frac{1}{s} \,
												{ |  \widehat{f} (s, \xi) | }^2  \widehat{f} (s, \xi) \, ds 
												+ \int_1^t R (s, \xi) \, ds
\end{align}
with
\begin{equation}
\label{R}
R (s,\xi) :=  - \frac{i}{s}  {(2\pi)}^{-1} 
\int  \left( e^{-i \eta \sigma/s } - 1 \right)  \F_{\eta, \sigma}^{-1} \left[  F \right] (s, \eta, \sigma, \xi) 
    										\, d\eta d\sigma
\end{equation}
We now claim that $R$ satisfies
\begin{equation}
\label{estimateR}
| R(s, \xi)  | \lesssim s^{- 1 - \d + 3\alpha} {\| u \|}^3_{X_T}
\end{equation}
for any $3 \alpha < \d < \frac{1}{2}$, uniformly in $\xi$.
We postpone for the moment the proof of this fact and continue the proof of the $L^\infty$
decay estimate.
Taking a derivative with respect to time in \eqref{intidfhat} one obtain
\begin{equation}
\label{partialtfhat}
\partial_t \widehat{f}(t,\xi) = - \frac{ i }{ t }
{ |  \widehat{f} (t, \xi) | }^2  \widehat{f} (t, \xi) +   R(t , \xi) \, .
\end{equation}
Since $R$ has a better decay in time, it is clear that one should consider
\begin{equation}
\label{what}
\widehat{w} (t,\xi) = \widehat{f} (t,\xi)  B (t,\xi) \hskip10pt \mbox{with} \hskip8pt  
										B (t,\xi)  = \exp \left( i \int_1^t \frac{{ |  \widehat{f} (s, \xi) | }^2}{s}  \, ds \right)
										\, ,
\end{equation}
so that
\begin{equation*}
\partial_t \widehat{w} (t,\xi) =  B(t,\xi) R(t,\xi)  \, ,
\end{equation*}
hence,
\begin{align*}
| \widehat{f} (t,\xi) |  =  | \widehat{w} (t,\xi) | 
												& \leq 
												| \widehat{u}_\ast (\xi) |  +  \int_1^t |  \partial_t \widehat{w} (s, \xi) | \, ds
												\\
												& \leq
												| \widehat{u}_\ast (\xi) |  +  \int_1^t |  R (s, \xi)  | \, ds \lesssim 
												{\| u_\ast \|}_{H^{0,1}}   +  t^{-\d + 3\alpha}  {\| u \|}^3_{X_T}  
												\, .
\end{align*}
Using lemma \ref{statNLS} with $n=1$, $3\alpha < \d$  and $\alpha < \beta < \frac{1}{4}$, one gets
\begin{equation*}
{\| u \|}_{L^\infty} \lesssim \frac{1}{t^{1/2}} \left( {\| u_\ast \|}_{H^{0,1}}   
													+  t^{-\d+3\alpha}  {\| u \|}^3_{X_T} \right)  +  \frac{1}{t^{{1/2} + \beta}}  {\| f \|}_{H^{0,1}}
													\lesssim
													\frac{1}{t^{1/2}} \left( {\| u_\ast \|}_{H^{0,1}}   +  {\| u \|}^3_{X_T} \right)   \, .
\end{equation*}
This concludes the proof of \eqref{cubicest} and implies
global existence for solutions of \eqref{eqs}--\eqref{gNLS}.

\vskip5pt
\nl
{\it Estimate for the remainder \eqref{estimateR}}:
From the definition of $R$ we immediately get
\begin{equation*}
| R (s,\xi) | \lesssim s^{-1-\d}  \int {|\sigma |}^\d {|\eta |}^\d   
																\left| \F^{-1}_{\eta, \sigma} [F] (s, \eta, \sigma, \xi)    \right|
																\, d\s d\eta
\end{equation*}
for any $0 < \d < \frac{1}{2}$.
To estimate this term we write an explicit expression for
$\F^{-1}_{\eta, \sigma} [F]$.
Since
\begin{eqnarray*}
  \mathcal{F}^{-1}_\eta \bigl[ \widehat{f}(\xi-\eta) \bigr] (\eta^\p) 
  = 
  e^{i \eta^\p \xi} f(-\eta^\p) \quad \mbox{and}
  \quad
  \mathcal{F}^{-1}_\eta \left[ \overline{\widehat{f}(\xi-\eta-\sigma)} \right] (\eta^\p)
  = 
  e^{i \eta^\p (\xi-\sigma)} \overline{f} (\eta^\p) \, ,
\end{eqnarray*}
one has
\begin{align*}
\mathcal{F}^{-1}_{\eta} \left[ \widehat{f}(\xi-\eta)\overline{\widehat{f}(\xi-\eta-\sigma)} \right] (\eta^\p)
  & =
  (2\pi)^{1/2} \mathcal{F}^{-1}_{\eta} [\widehat{f}(\xi-\eta)] \ast 
  \mathcal{F}^{-1}_{\eta} \left[  \overline{\widehat{f}(\xi-\eta-\sigma)}   \right] (\eta^\p)
	\\
  & = (2\pi)^{1/2} \int e^{i(\eta^\p - x) \xi} f(x - \eta^\p) \,
  e^{ix (\xi-\sigma)} \overline{f}(x) \,dx
  \\
  & = (2\pi)^{1/2} e^{i \xi \eta^\p}
    \int e^{-i x \sigma}f(x - \eta^\p) \overline{f}(x) \,dx
\end{align*}
Moreover, since
$
  \mathcal{F}^{-1}_{\sigma} [ e^{-ix \sigma}\,\widehat{f}(\xi-\sigma) ] (\sigma^\p)
  = e^{i\xi (\sigma^\p + x)} f(x - \sigma^\p)
$
we have
\begin{align}
\nonumber
  \mathcal{F}^{-1}_{\eta, \sigma} [F] 
  & = \mathcal{F}^{-1}_{\eta,\sigma} \left[
    \widehat{f}(\xi-\eta) 
    \overline{\widehat{f}(\xi-\sigma-\eta)} \widehat{f}(\xi-\sigma)
    \right]
    \\
    \nonumber
  & = (2\pi)^{1/2} e^{i\xi \eta^\p}\,
    \mathcal{F}^{-1}_{\sigma} \left[
    \int e^{-ix \sigma} f(x - \eta^\p) \overline{f}(x) \,dx \, \widehat{f}(\xi-\sigma) \right]
    \\
  & = (2\pi)^{1/2} e^{i\xi (\eta^\p + \sigma^\p) } 
  \int e^{-ix \xi}
  f(x - \eta^\p) \overline{f}(x) f(x-\sigma^\p) \,dx.
\label{explicitF}  
\end{align}
Therefore, 
\begin{eqnarray*}
  |R (s, \xi)| \lesssim 
  s^{-1-\delta} \int {|\eta|}^{\delta} \, {|\sigma|}^{\delta} \,
  |f(s,x - \eta)| |\overline{f}(s,x)| |f(s, x-\sigma)| \,dx \,d\eta \,d\sigma.
\end{eqnarray*}
%
%
%
%
Using $\d < \frac{1}{4}$ and
$
  {|\eta|}^{\delta} {|\s|}^{\delta}
  \lesssim \left( {|x - \eta|}^{\delta}+ {|x|}^{\delta} \right) \left( {|x - \s|}^{\delta}+ {|x|}^{\delta} \right)
$,
we can conclude that
\begin{eqnarray*}
  |R(s,\xi)| \lesssim  s^{- 1 - \delta} \, {\left\| {\langle x \rangle}^{2\delta} f \right\|}^2_{L^1} {\| f \|}_{L^1}
  \lesssim s^{- 1 - \delta}  \,  {\| f \|}_{H^{0,1}}^{3}
  \lesssim s^{ - 1 - \d + 3\alpha} {\| u \|}^3_{X_T} 
  \quad _\Box
\end{eqnarray*}

\section{Proof of proposition \ref{proH}}
\label{secproH}
{\it $\dot{H}^{m,0}$ estimates:}
For any ${|k|} = m$ and $p,q \geq 1$ with $\frac{1}{p} + \frac{1}{q} = \frac{1}{2}$ we have
\begin{align}
\label{Hder1}
{ \| \partial^k C \|}_{L^2} & \lesssim \int_1^t  { \left\|  {|x|}^{-1} \ast {|u|}^2 \right\|}_{L^\infty} 
																							 { \| \partial^k  u \|}_{L^2} \, ds
\\					
\label{Hder2}																	 
& + \int_1^t  { \left\|  {|x|}^{-1} \ast \partial^k {|u|}^2 \right\|}_{L^q} 
																							 { \| u \|}_{L^p} \, ds \, .
\end{align}
%
%
%
Using the endpoint of the Hardy--Littlewood--Sobolev inequality and \eqref{Lorentzdecay}, we get
\begin{align*}
\int_1^t  { \left\|  {|x|}^{-1} \ast {|u|}^2 \right\|}_{L^\infty} 
																							 { \| \partial^k  u \|}_{L^2} \, ds																							 & \lesssim 
\int_1^t  { \left\| u \right\|}^2_{L^{\frac{2n}{n-1},2}} 
																							 { \| u \|}_{X_T} s^{\alpha} \, ds
\\						
& \lesssim { \| u \|}^3_{X_T} \int_1^t  s^{-1} s^{\alpha} \, ds  \lesssim  { \| u \|}^3_{X_T} t^\alpha \, .
\end{align*}
To estimate \eqref{Hder2} we choose $q_1$ and $q_2$ satisfying $\frac{1}{q} = \frac{1}{q_1} - \frac{n-1}{n}$
and $\frac{1}{2} + \frac{1}{q_2} = \frac{1}{q_1}$, so that
\begin{eqnarray*}
{ \left\|  {|x|}^{-1} \ast \partial^k {|u|}^2 \right\|}_{L^q} \lesssim
{ \| \partial^k {|u|}^2 \|}_{L^{q_1}}  \lesssim 
{\| u \|}_{L^{q_2}} {\| \partial^k u \|}_{L^2} \, .
\end{eqnarray*}
Using again \eqref{Lorentzdecay} we obatin
\begin{align*}
\int_1^t  { \left\|  {|x|}^{-1} \ast \partial^k {|u|}^2 \right\|}_{L^q}  { \| u \|}_{L^p} \, ds 
& \lesssim
\int_1^t  {\| u \|}_{L^{q_2}} {\| \partial^k u \|}_{L^2}  { \| u \|}_{L^p} \, ds
\\
& \lesssim
{\| u \|}_{X_T}^3  \int_1^t  s^\alpha s^{ - n ( 1 - 1/p - 1/q_2 ) } \, ds \lesssim  {\| u \|}_{X_T}^3  t^\alpha \, .
\end{align*}

\vskip5pt
\nl
{\it $\dot{H}^{0,m}$ estimate in dimension $n \geq 3$:} 
Let $k$ be any multi--index with $|k| = m$. 
Applying $\partial_\xi^k$ to  $\widehat{C}$ produces a linear combination of terms of the form 
\begin{eqnarray}
\label{Hweight1}
& & \int_{1}^{t} \!\int 
    e^{is\eta\cdot\sigma} |\eta|^{-n+1}
    \partial_{\xi}^{k} \widehat{f} (s,\xi-\eta) \, 
    \left[ \overline{\widehat{f}(s,\xi-\eta-\sigma)} \widehat{f}(s,\xi-\sigma) \right] 
    \, d\eta d\sigma \, ds
\\
\label{Hweight2}
& & \int_{1}^{t} \!\int 
    e^{is\eta\cdot\sigma} |\eta|^{-n+1}
    \partial_{\xi}^{k- j} \widehat{f} (s,\xi-\eta) \, 
    \partial_{\xi}^{ j}  \left[ \overline{\widehat{f}(s,\xi-\eta-\sigma)} \widehat{f}(s,\xi-\sigma) \right] 
    \, d\eta d\sigma \, ds
\end{eqnarray}
for $0 \leq |j| \leq m$.
Since 
\begin{equation}
\label{byparts}
\partial_{\xi}^{j}  \left[ \overline{\widehat{f}(s,\xi-\eta-\sigma)} \widehat{f}(s,\xi-\sigma) \right] 
=
{(-1)}^{|j|} \partial_{\s}^{j}  \left[ \overline{\widehat{f}(s,\xi-\eta-\sigma)} \widehat{f}(s,\xi-\sigma) \right] 
\end{equation}
we can integrate by parts, with respect to $\s$, transforming \eqref{Hweight2} 
into a combination of terms of the form
\begin{equation}
\label{Hweight3}
\int_{1}^{t} \!\int 
    e^{is\eta\cdot\sigma} {(s\eta)}^{j} {|\eta|}^{-n+1}
    \partial_{\xi}^{k- j} \widehat{f} (s,\xi-\eta)\, 
    \overline{\widehat{f}(s,\xi-\eta-\sigma)}
    \widehat{f}(s,\xi-\sigma)\,d\eta d\sigma \, ds 
\end{equation}
Using again the Hardy--Littlewwod--Sobolev inequality, we can estimate
\begin{align*}
{\| \eqref{Hweight1} \|}_{L^2} & \lesssim 
			\int_{1}^{t}  {\left\| |x|^{-1} \ast |u(s)|^{2} \right\|}_{L^{\infty}}
      {\| e^{is\Delta/2} x^{k} f(s)\|}_{L^{2}} \, ds
\\
& \lesssim
      \int_{1}^{t}  {\|u(s)\|}_{L^{\frac{2n}{n-1},2}}^{2} \, s^{\alpha} {\| u \|}_{X_T} \,ds
\\
& \lesssim {\|u\|}_{X_T}^{3} \int_{1}^{t}   s^{-1} s^{\alpha} \,ds 
  \lesssim t^{\alpha}\, {\|u\|}_{X_T}^{3} \, .
\end{align*}
To estimate \eqref{Hweight3} we need the following:
\begin{lem}
\label{lemmaweight}
Let $1 \le j \le m$ and let $ 2 < p < 2n/(n-2j)$, then
\begin{equation*}
  {\| e^{is\Delta/2} |x|^{m-j}f \|}_{L^{p}}
  \lesssim  s^{-n (\frac{1}{2}-\frac{1}{p}) }  \left( {\|f\|}_{L^{2}} + {\| {|x|}^{m} f \|}_{L^{2}}  \right).
\end{equation*}
\end{lem}

\pr
Since $p>2$, we have
$$
  {\| e^{is\Delta/2} {|x|}^{m-j}f \|}_{L^{p}}
  \lesssim s^{-n(\frac{1}{2}-\frac{1}{p})}
    {\| {|x|}^{m-j} f \|}_{L^{p^\p}}.
$$
Using the upperbound on $p$, we can estimate
\begin{eqnarray*}
{\| {|x|}^{m-j} f \|}_{L^{p'}}
 		& \le & {\| \chi_{\{|x|\le 1\}} {|x|}^{m-j} f \|}_{L^{p^\p}}
    + {\| \chi_{\{|x|>1\}} {|x|}^{m-j} f \|}_{L^{p^\p}}
\\
  	& \le &  {\| \chi_{\{|x|\le 1\}} f \|}_{L^{p^\p}}
    + {\| \chi_{\{|x|>1\}} {|x|}^{-j} \|}_{L^{\frac{2p}{p-2}}}
      {\| {|x|}^m f \|}_{L^{2}}
\\
& \lesssim & {\| f \|}_{L^2} + {\| {|x|}^m f \|}_{L^2}   \,\, _\Box
\end{eqnarray*}
%
%
%
%
Now we choose $p,q$ and $q_1 \geq 1$ satisfying the following relations:
\begin{eqnarray*}
\frac{1}{p} + \frac{1}{q} = \frac{1}{2} 
\quad , \quad  
2 < p < \frac{2n}{n-2j}
\quad , \quad
\frac{1}{q_1} = \frac{1}{q} - \frac{| j| - n + 1}{n} \, .
\end{eqnarray*}
Then, from lemma \ref{lemmaweight} and Sobolev's embedding, we obtain
\begin{align*}
{\| \eqref{Hweight3} \|}_{L^2}  & \lesssim
			\int_{1}^{t} s^{|j|}  \, { \left\| \left( \partial^j  {|x|}^{-1} \ast {|u(s)|}^2 \right)
      e^{is\Delta/2} x^{k-j} f(s) \right\|}_{L^{2}}  \, ds
\\
& \lesssim
			\int_{1}^{t} s^{|j|} \, {\left\| {|\nabla|}^{|j| - n + 1} {|u(s)|}^{2} \right\|}_{L^{q}}
      {\|  e^{is\Delta/2}x^{k - j}  f(s)\|}_{L^{p}} \, ds \, .
\\
& \lesssim
    \int_{1}^{t} s^{| j|}\,
    {\| u(s) \|}^2_{L^{2q_1}}
    s^{-n(\frac{1}{2}-\frac{1}{p})}
    s^\alpha  {\| u \|}_{X_T} \,ds
\\
& \lesssim  {\|u\|}_{X_T}^{3}
    \int_{1}^{t} s^{| j|} \,
    s^{-n (1 - \frac{1}{q_1}) }
    s^{-n (\frac{1}{2}-\frac{1}{p})}
    s^{\alpha} \,ds
\\
& =  \int_1^t s^{-1} s^\alpha \, ds {\|u\|}_{X_T}^{3} \lesssim t^\alpha {\|u\|}_{X_T}^{3}
\end{align*}


\vskip5pt
\nl
{\it $\dot{H}^{{0,2}}$ estimate in dimension $n = 2$:}
For any $k \in \mathbb{N}^2$ with ${|k|}=2$, $\partial_\xi^k \widehat{C}$ is a combination of terms of the form
\begin{eqnarray}
\label{2dweight1}
& & \int_{1}^{t} \int 
    e^{is\eta\cdot\sigma} |\eta|^{-1} 
    \partial_{\xi}^2
      \widehat{f} (s,\xi-\eta)\,  
    \overline{\widehat{f}(s,\xi-\eta-\sigma)}
    \widehat{f}(s,\xi-\sigma)
    \,ds\, d\eta\, d\sigma 
\\
& &    
\label{2dweight2}
\int_{1}^{t} \int 
    e^{is\eta\cdot\sigma} |\eta|^{-1} 
    \partial_{\xi}
      \widehat{f} (s,\xi-\eta)\, 
    \partial_{\xi}
    \{\overline{\widehat{f}(s,\xi-\eta-\sigma)}
    \widehat{f}(s,\xi-\sigma)\}\,ds\, d\eta\, d\sigma 
\\    
\label{2dweight3}
& & \int_{1}^{t} \int 
    e^{is\eta\cdot\sigma} |\eta|^{-1} 
      \widehat{f} (s,\xi-\eta)\, 
    \partial_{\xi}^2
    \{\overline{\widehat{f}(s,\xi-\eta-\sigma)}
    \widehat{f}(s,\xi-\sigma)\}\,ds\, d\eta\, d\sigma  \, .
\end{eqnarray}
%
%
%
%
The term in \eqref{2dweight1} can be estimated in the same way as \eqref{Hweight1}, so we skip it.
Using \eqref{byparts} again,
we can integrate by parts transforming \eqref{2dweight2} and \eqref{2dweight3} in terms of the form
\begin{align}
\label{2dweight2.2}
&  \int_{1}^{t} \int 
    s  e^{is\eta\cdot\sigma} \frac{\eta}{|\eta|}
    \partial_{\xi}
      \widehat{f} (s,\xi-\eta)\, 
    \{\overline{\widehat{f}(s,\xi-\eta-\sigma)}
    \widehat{f}(s,\xi-\sigma)\} \, d\eta\, d\sigma \,ds
\\    
\label{2dweight3.2}
& \int_{1}^{t} \int 
    s  e^{is\eta\cdot\sigma} \frac{\eta}{|\eta|} 
      \widehat{f} (s,\xi-\eta)\, 
    \partial_{\xi}
    \overline{\widehat{f}(s,\xi-\eta-\sigma)}
    \widehat{f}(s,\xi-\sigma) \, d\eta\, d\sigma \,ds  
\\
\nonumber
& + \quad \mbox{similar term} \, .
\end{align}
Denoting by $\mathcal{R}$ the Riesz--transform, i.e.
$ \mathcal{R} f  := \F^{-1} ( \frac{\eta}{|\eta|} \F {f} )$, we have
\begin{align*}
{\| \eqref{2dweight2.2} \|}_{L^2} 
						& =  {\left\| \int_1^t s \left( \mathcal{R} {|u|}^2 \right) 
						\left( e^{-is \Delta/2} x f \right) \right\|}_{L^2} \, ds
						\\
						& \lesssim \int_1^t s \, {\| {|u|}^2 \|}_{L^4}  
						\, s^{-1/2} {\| x f \|}_{L^{4/3}} \, ds
						\\
						& \lesssim
						\int_1^t s^{-1} \, {\| u \|}_{X_T}^2  {\| {\langle x \rangle}^2 f \|}_{L^2} \, ds
						\lesssim 
						t^\alpha {\| u \|}_{X_T}^3 \, .
\end{align*}
For \eqref{2dweight3.2} we get 
\begin{align*}
{\| \eqref{2dweight3.2} \|}_{L^2} 
						& =  {\left\| \int_1^t s  \mathcal{R} 
						\left( \bar{u} e^{- is \Delta/2} (x f) \right) 
						e^{- is \Delta/2} f \right\|}_{L^2} \, ds
						\\
						& \lesssim \int_1^t s \, {\| \bar{u} e^{-is \Delta/2} (x f)  \|}_{L^2}  
						{\| e^{- is \Delta/2} f \|}_{L^\infty} \, ds
						\\
						& \lesssim
						\int_1^t   s^{-1} {\| u \|}_{X_T}^2  {\| x f \|}_{L^2} \, ds
						\lesssim 
						t^\alpha {\| u \|}_{X_T}^3 \, .
\end{align*}

\vskip5pt
\nl
{\it $L^\infty$ decay estimate}:
We proceed as in the NLS case applying Plancherel's identity to \eqref{int3}.
Since
\begin{align*}
  \mathcal{F}_{\eta} \left[  e^{is\eta\cdot\sigma}  {|\eta|}^{-n+1} \right] (\eta^\p)
  = \mathcal{F}_{\eta}  \left[{|\eta|}^{-n+1} \right] (\eta^\p - s\sigma)
  = C_{n-1} \, {|\eta^\p - s\sigma|}^{-1}
\end{align*}
we have
\begin{align*}
  \mathcal{F}_{\eta,\sigma}  \left[ e^{is\eta\cdot\sigma}  {|\eta|}^{-n+1} \right] (\eta^\p,\sigma^\p)
  & =   C_{n-1} \, \mathcal{F}_{\sigma}  \left[ |\eta^\p - s \sigma|^{-1} \right] (\sigma^\p)
  = s^{-1}\,e^{i \eta^\p \cdot \sigma^\p /s} 
    {|\sigma^\p|}^{-n+1} \, .
\end{align*}
Therefore, we can write
\begin{align*}
\widehat{C}(u,u,u) 
& = -i  {(2\pi)}^{-n/2} C_1 \int_1^t \int e^{is\eta\cdot\sigma}  {|\eta|}^{-n+1} F(s,\xi,\eta,\sigma)\,
    d\eta d\sigma 
\\
& = - i {(2\pi)}^{-n/2} C_1 \int_1^t \int \frac{1}{s}  e^{i \eta \cdot \sigma/s} 
    {|\sigma|}^{-n+1} \F_{\eta, \sigma}^{-1} [F](s,\xi,\eta,\sigma)\,
    d\eta d\sigma 
\\
& =: \int_1^t I_{0}(s, \xi)  +  \int_1^t R (s, \xi),
\end{align*}
with
\begin{align}
\nonumber
& I_0 (s, \xi) :=  - \frac{i}{s} {(2\pi)}^{-n/2} C_1 
		\int  {|\sigma|}^{-n+1} \F^{-1}_{\eta, \s}[F](s,\xi,\eta,\sigma)\,
    d\eta d\sigma  \, ,
\\
\label{RH}
& R (s,\xi) :=  - \frac{i}{s} {(2\pi)}^{-n/2} C_1    
		\int \left( e^{i\eta \cdot \sigma /s}  - 1 \right)  {|\s|}^{-n+1} 
		\F^{-1}_{\eta,\s} [F](s,\xi,\eta,\sigma)\,
    d\eta d\sigma
\end{align}
and $F$ given by \eqref{F}.
Next, we observe that
\begin{align*}
I_{0}(s,\xi) & = - \frac{i}{s} {(2\pi)}^{-n/2} C_1 \int \mathcal{F}_{\sigma}[ {|\sigma|}^{-n+1}] \, 
    \F_{\eta}^{-1} [F] (s,\xi,\eta,\sigma)
    \,d\eta\,d\sigma
\\
& =  - \frac{i}{s} {(2\pi)}^{-n/2}
		\int {|\sigma|}^{-1}
    \F_{\eta}^{-1} [F](s,\xi,\eta,\sigma)
    \,d\eta\,d\sigma
\\
& = - \frac{i}{s}
		\int {|\sigma|}^{-1}  F(s,\xi,0,\sigma)\,d\sigma
\\ 
& = - \frac{i}{s} 
		\int {|\sigma|}^{-1}
    {| \widehat{f}(s,\xi-\sigma) |}^{2}\,d\sigma\, \widehat{f}(s,\xi)
\\
& = - \frac{i}{s} 
		\left( {|x|}^{-1}\ast {|\widehat{f}(s)|}^{2} \right) (\xi) \, \widehat{f}(s,\xi) \, .
\end{align*}
Therefore,
\begin{eqnarray}
\label{intidfhatH}
\widehat{f} (t, \xi)  =  \widehat{u}_{\ast}(\xi)
    										- i  \int_1^t \frac{1}{ s }
												\left( {|x|}^{-1}\ast {|\widehat{f}(s)|}^{2} \right) (\xi)  \widehat{f} (s, \xi) \, ds 
												+ \int_1^t R (s, \xi) \, ds
\end{eqnarray}
with $R$ given by \eqref{RH}.
We now claim that $R$ satisfies
\begin{equation}
\label{estimateRH}
| R(s, \xi)  | \leq s^{- 1 - \d + 3\alpha} {\| u \|}^3_{X_T}
\end{equation}
for any $3 \alpha < \d < \frac{1}{2}$, $\xi \in \R^n$. 
Assuming for the moment \eqref{estimateRH},
we take a derivative with respect to time in \eqref{intidfhat} and obtain
\begin{equation}
\label{partialtfhatH}
\partial_t \widehat{f}(t,\xi) = - \frac{ i }{ t }
										\left( |x|^{-1}\ast |\widehat{f}(t)|^{2} \right)  \widehat{f} (t, \xi) +   R(t , \xi) \, .
\end{equation}
Since $R$ has a better decay in time we eliminate the middle term in \eqref{partialtfhat} 
by considering a modified profile
\begin{equation}
\label{whatH}
\widehat{w} (t,\xi) = \widehat{f} (t,\xi)  B (t,\xi) \hskip10pt \mbox{with} \hskip8pt  
										B (t,\xi)  = 
										\exp \left( i \int_1^t  |x|^{-1}\ast |\widehat{f}(s)|^{2}  \, ds \right) \, .
\end{equation}
Then
$\partial_t \widehat{w} (t,\xi) =  B(t,\xi) R(t,\xi)$,
whence
\begin{eqnarray*}
| \widehat{f} (t,\xi) |  =  | \widehat{w} (t,\xi) | \leq 
												| \widehat{u}_\ast (\xi) |   +  \int_1^t |  R (s, \xi)  | \, ds \lesssim 
												{\| u_\ast \|}_{H^{0,m}}   +  t^{-\d+3\alpha}  {\| u \|}^3_{X_T}  
												\, .
\end{eqnarray*}
Using lemma \ref{statNLS}, we finally obtain
\begin{equation*}
{\| u \|}_{L^\infty} \lesssim \frac{1}{t^{n/2}} \left( {\| u_\ast \|}_{H^{0,m}}   
													+  t^{-\d+3\alpha}  {\| u \|}^3_{X_T} \right)  
													+  \frac{1}{t^{{n/2} + \beta}}  {\| f \|}_{H^{0,m}}
													\lesssim
													\frac{1}{t^{n/2}} \left( {\| u_\ast \|}_{H^{0,m}}   +  {\| u \|}^3_{X_T} \right)   \, .
\end{equation*}
This concludes the proof of global existence of solutions 
of \eqref{eqs}--\eqref{gH}, provided we show \eqref{estimateRH}.


\vskip5pt
\nl
{\it Remainder estimate \eqref{estimateRH}}:
From the definition of $R$ we immediately get
\begin{equation*}
| R (s,\xi) | \lesssim s^{-1-\d}  \int {|\eta |}^\d {|\sigma |}^{-n + 1 + \d}
																\left| \F_{\eta, \sigma}^{-1} [F] (s, \eta, \sigma, \xi)   \right|
																\, d\s d\eta
\end{equation*}
for any $0 < \d < \frac{1}{2}$.
Now we use the explicit expression for
$\F_{\eta, \s}^{-1} [F] (s,\xi,\eta,\sigma)$, derived in \eqref{explicitF}, to obtain
\begin{eqnarray*}
  |R (s, \xi)| \lesssim 
  s^{-1-\delta} \int {|\eta|}^{\d} \, {|\sigma|}^{-n + 1 + \d}  \,
  |f(x-\eta)| |\overline{f}(x)| |f(x-\sigma)| \, dx \, d\eta \, d\sigma \, .
\end{eqnarray*}
Since $\d < \frac{1}{2}$, the integral with respect to $\sigma$ can be directly estimated by
\begin{eqnarray*}
  \int {|\sigma|}^{-n + 1 + \d}    |f(x-\sigma)  |  \, d \sigma
  \lesssim {\| f \|}_{L^\infty} + {\| f \|}_{H^{0,m}}  \lesssim {\| u \|}_{X_T} s^\alpha \,  .
\end{eqnarray*}
Using 
$
{|\eta|}^{\delta} \lesssim {|\eta+x|}^{\delta} + {|x|}^{\delta}
$
we conclude
\begin{eqnarray*}
  |R(s,\xi)| \lesssim  s^{- 1 - \delta + \alpha} \, {\| |x|^{\delta} f \|}_{L^1}  {\|f\|}_{L^{1}}  {\| u \|}_{X_T} 
  \lesssim s^{- 1 - \delta + 3\alpha}  \,  {\| u \|}_{X_T}^{3}
  \quad _\Box
\end{eqnarray*}

\addcontentsline{toc}{section}{Bibliography}


\end{document}